# NORMAL SURFACES AND COLORED KHOVANOV HOMOLOGY

CHRISTINE RUEY SHAN LEE

ABSTRACT. We show that colored Khovanov homology detects classes of essential surfaces as a direct analogue of the slope conjectures for the colored Jones polynomial. We do this by identifying certain generators of the colored Khovanov chain complex with normal surfaces in the complement of the knot using an ideal triangulation from a diagram.

## CONTENTS

## 1. INTRODUCTION

We study the colored Jones polynomial, a generalization of the Jones polynomial, and its categorification, colored Khovanov homology. Fix a complex number $q$ that is not a root of unity. To a knot $K \subset S^3$, the colored Jones polynomial assigns a sequence of Laurent polynomials $\{J_K^n(q)\}$ in $\mathbb{Z}[q, q^{-1}]$ indexed by natural numbers $n \geq 2$, where $J_K^2(q)$ is the Jones polynomial. In a different direction, based on the state sum model of the colored Jones polynomial, categorification assigns a bi-graded chain complex $\{CKh_{i,j}^n\}$ from which the $n$th colored Jones polynomial may be recovered from a suitable Euler characteristic of the homology groups.

A goal of quantum topology is to understand how quantum invariants such as the Jones polynomial and its categorification, Khovanov homology, encode geometric and topological information used to define other link invariants. Examples are the volume of hyperbolic knots [Kas97, MM01,

MMO$^{+}$02] and the set of boundary slopes [Gar11b, KT15]. A better understanding of this relationship may help construct new invariants of manifolds and provide new insights into quantum invariants.

In this paper, we propose a model for understanding the colored Jones polynomial and colored Khovanov homology via normal surfaces in the exterior of a knot. Our approach uses normal surface theory applied to a fixed diagrammatic triangulation of the knot exterior. We establish a relationship between certain terms in a state sum defining the polynomial and its categorifying chain complex, which we call "colored surface states," see Definition 3.17, and normal surfaces in the triangulation.

**Theorem 1.1.** *Let $K$ be a nontrivial knot in $S^3$ with diagram $D$. Suppose a colored Kauffman state $\sigma$ is a colored surface state with homological grading $h_\sigma$, then there is a corresponding normal surface $\mathcal{N}_\sigma$ in the octahedral triangulation $\mathcal{T}_D$ of the knot $K$ with slope $s_\sigma$ such that*

$$h_\sigma = s_\sigma n^2.$$

Normal surface theory plays a pivotal role in 3-dimensional topology where it has important applications [JS03a, HLP99] by studying the set of normal surfaces in a triangulation of a 3-manifold. In particular, the theory solves the problem of finding the set of boundary slopes of essential surfaces by enumerating fundamental normal surfaces and their slopes, since it is known that every essential surface is isotopic to a normal surface, and every normal surface can be obtained as a Haken sum of fundamental normal surfaces.

On the other hand, recent research ([GL05], [Gar11b], [Gar11a], [FKP11], [FKP13], [KT15], [GvdV16], [LvdV16], [MT17], [BMT18], [Lee], [LYL19], [GLvdV20]), beginning with the slope conjecture stated by Garoufalidis [Gar11b], suggests a close relationship between essential surfaces of a knot and the asymptotics of the degree of the colored Jones polynomial. For the knots studied, the essence of the relationship is captured by state-sum generators of colored Khovanov homology and the correspondence of their homological and quantum gradings to the boundary slopes of matching essential surfaces. In Theorem 1.1, we give a general criterion for matching a generator to normal surfaces which are candidates for essential surfaces.

These normal surfaces can have multiple sheets and non-integral slopes and thus form a strictly larger class than *state surfaces*, which are spanning surfaces from a Kauffman state for which the correspondence to generators in Khovanov homology is well known, see [Oza11, FKP11, Kin18]. Our second result generalizes these results to colored Khovanov homology by considering the homology class of a generator $X_\sigma$ from a colored surface state in colored Khovanov homology. In the theorem that follows, a knot diagram $D$ is *highly twisted* if every twist region of $D$ has more than two crossings.

**Theorem 1.2.** *Let $K$ be a nontrivial knot in $S^3$ with diagram $D$. Suppose that $D$ is highly twisted, and the normal surface $\mathcal{N}_\sigma$ corresponding to a colored surface state $\sigma$ is essential. Then the element $X_\sigma$ assigned to $\sigma$ is a cycle and its homology class $[X_\sigma]$ is nontrivial in colored Khovanov homology.*

Theorem 1.2 shows that colored Khovanov homology detects essential surfaces from colored surface states for highly twisted links which can be compared with results in [AL19] that show Bordered Floer homology detects incompressible surfaces; see also [Kin18], which says that homogeneously-adequate Kauffman states corresponding to essential surfaces give non-trivial Khovanov homology classes. The advantage of considering the colored Khovanov homology over the colored Jones polynomial is that more boundary slopes may be detected, not just the ones that realize the extremal quadratic growth rates of the degree of the colored Jones polynomial from the viewpoint of the slope conjecture [Gar11b, KT15]. However, to give a complete topological model, one should assign meaning to every generator in the colored Khovanov complex. To this end, we are not completely successful.

Our results suggest a new criterion, using quantum invariants, for finding essential surfaces of a knot using a diagram. This criterion applies to knots that are not Montesinos knots, for whom Hatcher and Oertel's algorithm [HO89] gives a complete classification of boundary slopes. In particular, we give an example of a knot that is not a Montesinos knot, where computational evidence by way of SnapPy [CDGW] shows that a normal surface $\mathcal{N}_\sigma$ detected by Theorem 1.1 indeed gives a boundary slope. We remark that Theorem 1.1 suggests a systematic approach to finding candidate normal surfaces for essential surfaces detected by colored Khovanov homology, and we will pursue the question of whether such normal surfaces are essential in a future project.

On this note, it also seems possible to directly compare the construction in this paper to that of the 3D-index [GHHR16, FKB08], which is an invariant defined on the set of ideal triangulations of a 3-manifold using the corresponding normal surfaces for the ideal triangulation. We would also like to mention that there are many models for Khovanov homology, see [Law90, Ang19, AS19, SS06], each with specific target invariants such as Lagrangian Floer homology [Ang19], symplectic Floer homology [AS19, SS06], and numerous others whose details are beyond the scope of this paper. Distinct from these approaches, our perspective presents yet another way to relate quantum invariants to 3-manifold topology.

**Organization.** The paper is organized as follows. In Section 2 we define the diagrammatic triangulation we use and related preliminary notions in normal surface theory. In Section 3 we describe the setup for colored Khovanov homology. Then we prove Theorem 1.1 in Section 3.7 and Theorem 1.2 in Section 3.8. In Section 4 we give an example of a non-Monteinos knot for which we can pinpoint a boundary slope suggested by Theorem 1.1.

## 2. A DIAGRAMMATIC TRIANGULATION OF THE KNOT EXTERIOR

Let $K \subset S^3$ be a knot and fix a diagram $D$ of $K$. We describe the diagrammatic triangulation $\mathcal{T}$ we use for $S^3 \setminus K$. The triangulation $\mathcal{T}$ is obtained by applying the inflation procedure in [JR14] to the diagrammatic ideal triangulation described in detail by [KKY18]. Then, we prove Lemmas 2.5, 2.6, and 2.8 giving a diagrammatic criterion for normal surfaces coming from normal disks in twist regions. This will be used to connect generators of colored Khovanov homology to normal surfaces in $S^3 \setminus K$ in Section 3.

An *ideal triangulation* is a triangulation of a 3-manifold by ideal tetrahedra. More precisely, let $\triangle = \{\triangle_1^3, \triangle_2^3, \ldots, \triangle_m^3\}$ be a disjoint collection of 3-simplices. Let $\Phi$ be a collection of face-pairing Euclidean isometries between 2-simplices $\triangle_i^2, \triangle_j^2$ of the 3-simplices of $\triangle$, where it is permissible that $i = j$. Consider the restriction $\Phi_0$ of $\Phi$ to $\triangle \setminus \triangle^0$, the collection of 3-simplices $\triangle$ with the 0-skeleton $\triangle^0$ removed. Let $int(M)$ be the interior of a compact 3-manifold $M$ with nonempty boundary, then $M$ is said to admit an *ideal triangulation* if $int(M)$ is homeomorphic to $(\triangle \setminus \triangle^0)/\Phi_0$ for some $\triangle$ and $\Phi$. The manifold $M$ is said to admit a *triangulation* if it is homeomorphic to $\triangle/\Phi$ for some collections of $\triangle$ and $\Phi$.

### 2.1. Octahedral ideal triangulation.
The ideal triangulation $\mathcal{T}^*$ we use for the exterior of a knot $K \subset S^3$ comes from a diagram $D$ of the knot as first described by Thurston in unpublished work. We will follow the detailed description of [KKY18]. See also [Wee05] for additional references on this ideal triangulation.

Position the knot in $S^2 \times I \subset S^3$ via a projection to the diagram $D$ in $S^2 \times 1/2$. Call the 3-balls in $S^3$ with boundary $S^2 \times \{0, 1\}$ the North pole $P_n$ and the South pole $P_s$, respectively. We will start with an ideal triangulation of the manifold $S^3 \setminus \{K \sqcup P_s \sqcup P_n\}$.

An ideal octahedron is an octahedron with vertices removed. We give it an ideal triangulation by subdividing it in 5 ideal tetrahedra, with gluing instructions given in Tables 1, 2. Note here $z$ is the center tetrahedron, $u_f$ is the upper front, $u_b$ is the upper back, $\ell_f$ is the lower front, and $\ell_b$ is the lower back tetrahedron.

| Tetrahedron | Face 012 | Face 013 | Face 023 | Face 123 |
|---|---|---|---|---|
| $z = o(1234)$ | $u_f(032)$ | $\ell_f(312)$ | $u_b(023)$ | $\ell_b(021)$ |
| $u_f = o(0123)$ | $u_b(012)$ | - | $z(021)$ - | - |
| $u_b = o(0134)$ | $u_f(012)$ | - $z(023)$ | - | - |
| $\ell_f = o(1245)$ | $\ell_b(301)$ | - | - | $z(130)$ |
| $\ell_b = o(2345)$ | $z(132)$ | $\ell_f(120)$ | - | - |

Table 1. Gluing instruction for an octahedron at a negative crossing ⟨✕⟩.

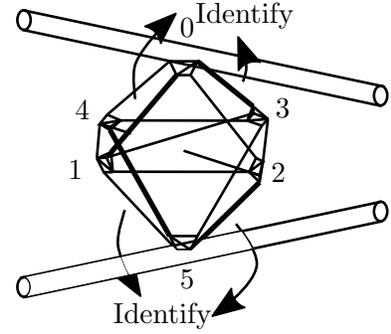

Figure 1

| Tetrahedron | Face 012 | Face 013 | Face 023 | Face 123 |
|---|---|---|---|---|
| $z = o(1234)$ | $u_f(231)$ | $\ell_f(320)$ | $u_b(102)$ | $\ell_b(230)$ |
| $u_f = o(0124)$ | $u_b(301)$ | - | - | $z(201)$ |
| $u_b = o(0234)$ | $z(203)$ | $u_f(120)$ | - | - |
| $\ell_f = o(1235)$ | $\ell_b(012)$ | - | $z(310)$ | - |
| $\ell_b = o(1345)$ | $\ell_f(012)$ | - | $z(312)$ | - |

Table 2. Gluing instruction for an octahedron at a positive crossing ⟨✕⟩.

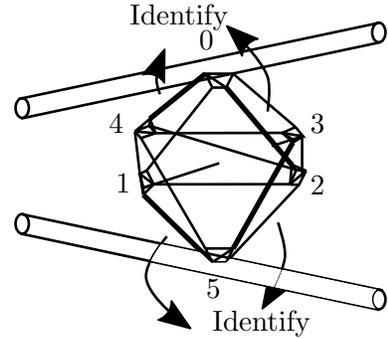

Figure 2

We specify the remaining isometries between faces of the octahedron $o$ for each negative/positive crossing as follows: For a negative crossing and the associated octahedron $o$, identify edges $o(01) \leftrightarrow o(03)$ and $o(25) \leftrightarrow o(45)$, so that a pair wraps around the top strand and another pair wraps around the bottom strand as indicated in Figure 1. Similarly for a positive crossing and the associated octahedron $o$ we identify edges $o(02) \leftrightarrow o(04)$ and $o(15) \leftrightarrow o(35)$ as indicated in Figure 2. Between a crossing and another joined by a strand in the link diagram, we identify the corresponding pair of faces by sliding them along the strand to meet, see Figure 3 for an example merging the faces of the octahedra of a negative crossing and a positive crossing. After we make these identifications connecting the octahedron at every crossing, we get an ideal triangulation $\mathcal{T}^\circ$ of $S^3 \setminus \{K \sqcup P_n \sqcup P_s\}$.

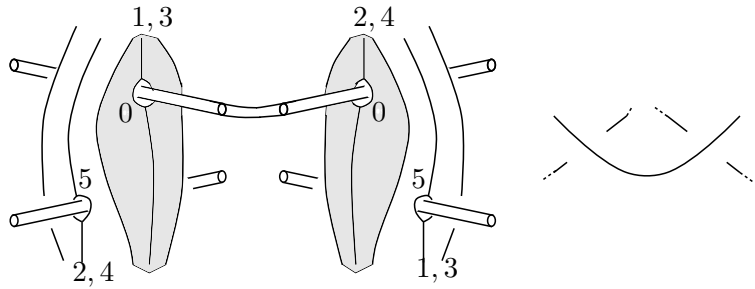

Figure 3. Joining faces (shaded grey) of two octahedra.

To obtain an ideal triangulation $\mathcal{T}^*$ of the knot exterior $S^3 \setminus K$, we add two ideally-triangulated "pillows" with a pre-drilled tube (gluing instructions in Table 3) to cancel the North pole and the South pole as described in [Wee05].

Pick two distinct 2-simplices $\triangle_j^2$, $\triangle_k^2$ in $\mathcal{T}^\circ$, so that $\triangle_j^2$ has a vertex on $P_n$ and another vertex on $N(K)$, and $\triangle_k^2$ has a vertex on $P_s$ and another vertex on $N(K)$. Let $P_{1,2}$ be the two triangular pillows. Discard the original face-pairing of one of the 2-simplices, say $\triangle_j^2 \leftrightarrow \triangle_{j'}^2$ and add in the triangular pillows in order by identifying one face of the triangular pillow $P_1$ to $\triangle_j^2$ and the other face of $P_1$ to $\triangle_{j'}^2$. Similarly we identify the faces of $P_2$ to $\triangle_k^2$ and $\triangle_{k'}^2$. This is oriented so that the two ends of the tube of the pillow match up with the two vertices, one on one of the poles and another one on $N(K)$. The interior of the resulting manifold with ideal triangulation $\mathcal{T}^*$ is then homeomorphic to $S^3 \setminus K$.

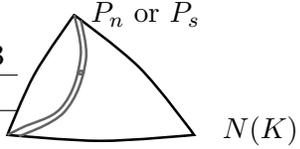

| Tetrahedron | Face 012 | Face 013 | Face 013 | Face 123 |
|---|---|---|---|---|
| 0 | 1 (013) | - | - | 1(012) |
| 1 | 0 (123) | 0 (012) | 1(123) | 1(023) |

TABLE 3. Gluing instructions for a pillow $P_{1,2}$ with a pre-drilled tube [HMP20].

FIGURE 4. Pillow.

## 2.2. A planar graph $G$ from a link diagram.
We will organize the tetrahedra in $\mathcal{T}^*$ and the tetrahedra added in the inflation procedure in the next section for the triangulation $\mathcal{T}$ using a graph $G$ that captures the twist region information of the knot diagram $D$.

Let $G$ be a planar, weighted, finite, and connected graph. Every knot or link $K$ admits a diagram $D$ such that $D = \partial F_G$, where $F_G$ is a surface obtained by replacing each vertex of $G$ by a disk in the plane, so that set of disks is disjoint, and by replacing each edge of $G$ by a twisted band, where we have right-handed/left-handed crossings depending on whether the weight on the edge is positive/negative, respectively. See Figure 5 for an example of how to recover the graph $G$ from a link diagram $D$, in which case we write $G = G(D)$, and the convention on the signs of the crossings.

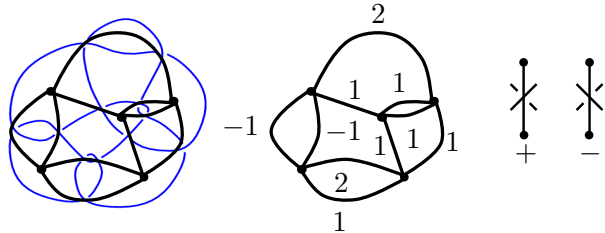

FIGURE 5. Weighted graph $G$ for a link diagram $D = \partial(F_G)$.

A *twist region* is a tangle diagram $T$ obtained from joining one or more crossings of the same sign. We consider twist regions in a link diagram corresponding to edges of $G$. Let $c(T)$ be the number of crossings in the twist region, we will label the crossings in a twist region corresponding to an edge $e$ in $G$ following the direction on $e$ as $c_1, c_2, \ldots, c_j, \ldots, c_{c(T)}$, and the respective octahedron $o_1, o_2, \ldots, o_j, \ldots, o_{c(T)}$.

## 2.3. Inflation of ideal triangulations.
From any ideal triangulation $\mathcal{T}^*$ of the interior of a compact 3-manifold $M$ with boundary, no components of which is a 2-sphere, Jaco and Rubinstein [JR14] showed one can obtain a triangulation $\mathcal{T}$ of $M$ by "inflating" the ideal triangulation $\mathcal{T}^*$. In this paper, we apply the procedure to the case where $M$ is the knot exterior $S^3 \setminus K$ with the octaheral ideal triangulation $\mathcal{T}^*$ discussed in Section 2.1, and we follow similar notation as in [JR14].

**Step 1**. Choose an index 4 *frame* $\Lambda$ in the torus boundary $\partial M$.

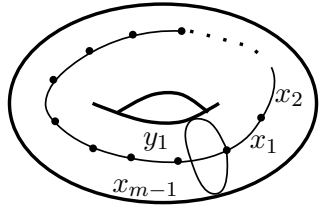

FIGURE 6. The index 4 frame on the torus.

The ideal triangulation $\mathcal{T}^*$ on $M$ induces a triangulation $\mathcal{T}^*_{\partial M}$ on $\partial M$. A graph in the 1-skeleton of the triangulation on $\partial M$ is a *spine* if each component of its complement in the torus is an open disk. A spine is a *frame* if it is minimal with respect to set inclusion. That is, if $\xi, \xi'$ are spines for $\mathcal{T}^*_{\partial M}$ and $\xi$ is a frame, then $\xi' \subset \xi$ implies $\xi = \xi'$.

The *index* of a vertex of a frame is the number of edges of the frame that meet at that vertex. If the index of a vertex is $> 2$, then the vertex is called a *branch point*. A *branch* of the frame is the closure of a component of the frame with branch points removed. Up to graph isomorphism, there are only two types of frames for the torus $T$: an index 4 frame or a double index 3 frame. We will choose an index 4 frame. See Figure 6.

We assemble the frame by picking, for each square of $\mathcal{T}^*_{\partial M}$ that comes from a strand through a twist region, an edge in the 1-skeleton as indicated in Figure 7. This gives a natural direction on the branch from an orientation on the knot diagram, as well as a transverse direction selected by the right-hand rule (with the thumb pointing in the direction of the frame). We choose the frame so that between two crossings in a twist region, both branches consist of edges of the outer faces $o_j(034)$ and $o_j(125)$, and the edges of a frame between two twist regions consist of edges of the inner faces $o_j(125)$ and $o_j(145)$.

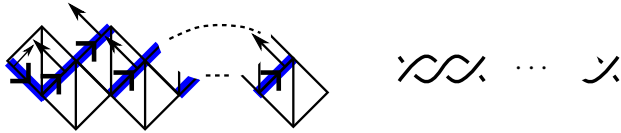

FIGURE 7. The choice of the frame (left) in blue and thicker highlighting with transverse direction is shown for a twist region (right).

After specifying the parts of the frame for each twist region, we pick a twist region and a square in the corresponding triangulation to place the branch point. The branch point is placed to avoid any intersection with the triangular pillows from the construction of $\mathcal{T}^*$ at the end of Section 2.1. This is always possible when the knot diagram has more than two crossings, as it is assumed throughout the paper that the knot is nontrivial. See Figure 8 below.

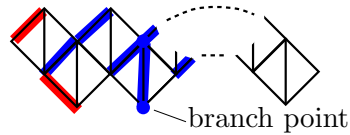

FIGURE 8. Edges of the frame are in blue and thicker highlighting and edges highlighted in red are where the pillows are inserted.

Label the edges of our chosen index 4 frame successively starting from the initial branch point. The branch which defines a longitude of the boundary torus is labeled $x$ with edges $x_1, \ldots, x_m$, and the branch which defines a meridian is labeled $y_1$ with a single edge.

**Step 2**. Inflate at every face of the triangulation $\mathcal{T}^*$ that intersects an edge of the frame.

For each face of $\mathcal{T}^*$ that intersects an edge of the frame, it intersects in 1, 2, or 3 edges. In our case, every face of $\mathcal{T}^*$ intersects the frame in a single edge, since by construction in Step 1, there are no intersections of the frame except at vertices with the North pole, the South pole, or the inserted pillows, where a face could intersect in 2 or more edges.

Inflation at an edge of the frame means we add a new tetrahedron to what will become the set of tetrahedra for the new triangulation $\mathcal{T}$. Let $(p)(abc) \in \mathcal{T}^*$ be a face with face pairing $(p)(abc) \leftrightarrow (p')(a'b'c')$ that intersects the frame in a single edge, say $x_j$. We discard the face pairing $(p)(abc) \leftrightarrow (p')(a'b'c')$ and add a new tetrahedron denoted $(x_j)$ corresponding to the edge $x_j$ of the frame with vertices $0, 1, 2, 3$. Then we make new face pairings depending on whether the transverse direction to the frame points into or out of $(p)$. If the transverse direction to the frame points out of $(p)$, then the new face pairings are $(p)(abc) \leftrightarrow (x_j)(032)$ and $(x_j)(132) \leftrightarrow (p')(a'b'c')$. Otherwise, if the transverse direction to the frame points into $(p)$, then the new face pairings are $(p)(abc) \leftrightarrow (x_j)(132)$ and $(x_j)(032) \leftrightarrow (p')(a'b'c')$. See Figures 9 and 10 for an illustration of these two cases.

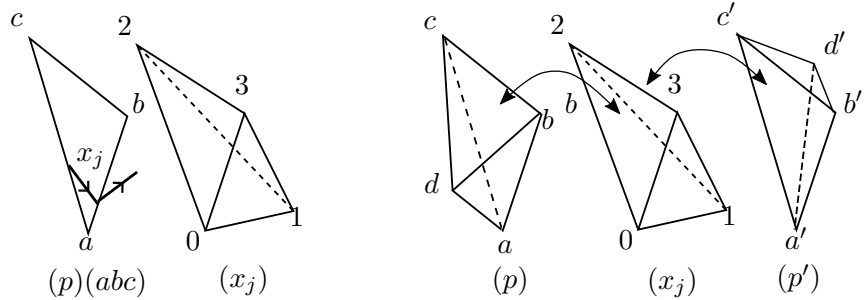

FIGURE 9. Identifying two faces of the newly added tetrahedron through inflating at a face. The case shown has the transverse orientation pointing out of $(p)$.

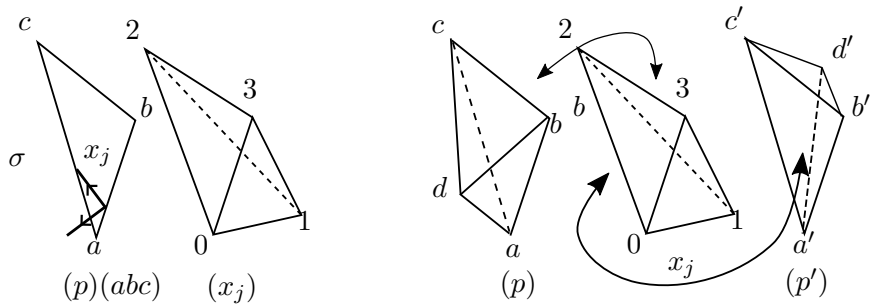

FIGURE 10. Identifying two faces of the newly added tetrahedron through inflating at a face. The case shown has the transverse orientation pointing into $(p)$.

We thus inflate at every edge of the frame. In our case with the index 4 frame of the torus, we add a new tetrahedron for every edge $x_1, \ldots, x_m$ and $y_1$. Note that for each of these new tetrahedron, two faces are unidentified with another face.

**Step 3**. Identify the faces of the newly-created tetrahedra from Step 2 by inflating at an edge.

Recall that from Step 2, there are unidentified faces of the newly added tetrahedra. In this step we make the remaining face-pairings.

We consider an edge, say $e \in \mathcal{T}^*$, whose endpoints $e_-$ or $e_+$ intersect the frame at a vertex. Let $D_e^-$ and $D_e^+$ be small regular neighborhoods of $e_+$ and $e_-$ in the boundary torus of $M$. There is a naturally cyclic order of the edges of the frame with vertex at $e_-$ or $e_+$. In [JR14], Jaco and Rubinstein names 3 configurations for the pattern of edges of the frame at one such vertex called *generic*, *crossing*, and *branch*. They show that every possible configuration can be reduced to a combination of these three, and how to identify the newly-created faces from Step 2 for those configurations.

For our chosen frame $\Lambda$, there are two edges of $\mathcal{T}^*$ which each intersects the index 4 branch point at a single point, corresponding to the branching configuration. The rest of the edges of $\mathcal{T}^*$ either does not intersect the frame, or only intersects the frame at a single point of index 2, corresponding to the generic configuration.

**Case 1: The generic configuration**. In the generic configuration the edge $e$ meets the frame in a single point that is a vertex of index 2 in the frame. Suppose without loss of generality that $e$ meets the frame in $e^-$, which is the vertex between the edges $x_j$ and $x_{j+1}$ of the frame. We make the new face pairing between the previously identified faces of $(x_j)$ and $(x_{j+1})$: $(x_j)(013) \leftrightarrow (x_{j+1})(012)$.

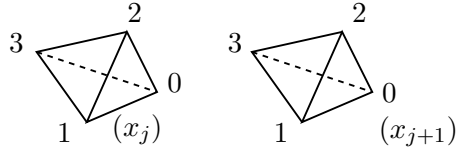

Figure 11. The generic configuration.

**Case 2: The branching configuration**. In the branching configuration the edge $e \in \mathcal{T}^*$ meets the frame in one point that is a branch point of index $b \geq 3$. Suppose again without loss of generality that $e$ meets the frame at $e^-$. In $D_e^-$ there are $b$ edges in the frame, ordered cyclically as $x, y, \ldots, z$. There are $b$ unidentified faces of the tetrahedra about the edge $e$ from inflating faces intersecting edges of each branch. Let $P_b$ be a planar $b$-gon and form the cone $0 \star P_b$ over the $b$-gon $P_b$. The cone $0 \star P_b = b^*$ has $b$ triangular faces: $(b^*)(120), (b^*)(230), \ldots, (b^*)(b10)$. We make the face identifications $(x_1)(012) \leftrightarrow (b^*)(b10)$, $(y_k)(013) \leftrightarrow (b^*)(210), \ldots$.

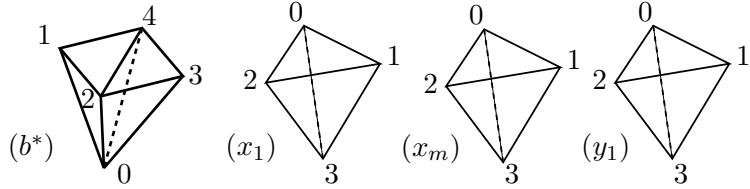

Figure 12. The branching configuration. In our case there is a single edge in the branch labeled $y_1$.

In our case, our branch point has index 4 so we create the planar 4-gon, and we make an arbitrary choice of the subdivision of the planar 4-gon $b^*$ as shown in Figure 12. The two remaining free faces $(b^*)(124)$ and $(b^*)(234)$ of $b^*$ form the boundary of the manifold. The face-pairings in $\mathcal{T}^*$ not involved in Steps 1-3 remain, and together with the newly-created tetrahedra and face-pairings they give a triangulation $\mathcal{T} = \mathcal{T}_\Lambda = (\triangle_\mathcal{T}, \Phi)$ of the manifold $M = S^3 \setminus K$ as stated below.

**Theorem 2.1.** [JR14] *Suppose $X$ is a compact 3-manifold with boundary, no component of which is a 2-sphere, and $\mathcal{T}^*$ is an ideal triangulation of the interior of $X$. Then for any collection of frames $\Lambda$, one frame in each of the vertex-linking surfaces of $\mathcal{T}^*$, the underlying point set of $\mathcal{T}$ is a compact 3-manifold $M$ homeomorphic to $X$ and the triangulation $T_\Lambda$ is a triangulation of $M$.*

2.4. **Normal surfaces and $Q$-matching equations.** With the triangulation $\mathcal{T}$ constructed in the previous section at hand, we proceed to construct normal surfaces in the triangulation. In this section we will collect fundamental results on normal surface theory and prove results used in the later sections about normal surfaces in our chosen triangulation of $S^3 \setminus K$.

A *normal disk* is a properly embedded disk in a 3-simplex $\triangle^3$ whose boundary is made up of properly embedded arcs in the faces of the 3-simplex, such that no two arcs are on the same face, and no two endpoints of any arc are on the same edge.

Up to normal isotopy of the 3-simplex, i.e., an isotopy of the 3-simplex which preserves all of its sub-simplices, there are exactly seven normal isotopy classes of normal discs in a 3-simplex. Four are normal triangular disks, each of which cuts off a corner of the tetrahedron, and three are normal quadrilateral disks, see Figure 13. We will call the image of an edge $e$ in a 3-simplex in $\mathcal{T}$ a *singlex*, and we will use the same notation for a simplex and its singlex, where it will be clear from the context whether we are referring to the image in the triangulation.

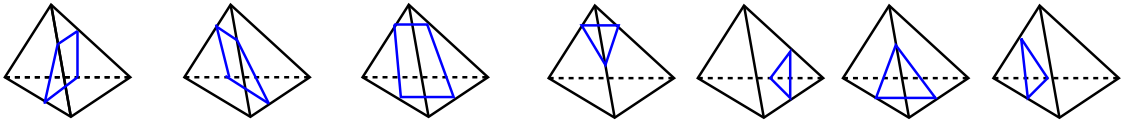

FIGURE 13. Left: Quadrilateral disks. Right: Triangular disks.

Let $M$ be a compact 3-manifold with a triangulation $\mathcal{T}$.

**Definition 2.2.** A *normal surface* of $M$ with respect to the triangulation $\mathcal{T}$ on $M$ is a subset of $M$ which intersects the 3-singlices in $\mathcal{T}$ in a pairwise disjoint collection of normal triangular and quadrilateral disks.

For a compact 3-manifold $M$ with triangulation $\mathcal{T}$ and $t$ tetrahedra, fix an ordering of the 3-singlices $\triangle_1, \ldots, \triangle_t$. The normal coordinates of a collection of normal disks is a $7t$-tuple:

$$(1) \qquad \vec{F} = (q_1', q_1'', q_1''', a_1', a_1'', a_1''', a_1'''', q_2', \ldots, a_t'''')$$

with nonnegative integer entries specifying the number of types of normal disks ($q$ for quadrilaterals and $a$ for triangular) in each 3-singlex. For example, $(0, 1, 0, 1, 0, 0, 0)$ indicates a quadrilateral disk of type 2 and a triangular disk of type 1.

For a 1-singlex $e$ in $\mathcal{T}$, consider the disjoint union $\tilde{B}(e)$ of all 3-singlices whose image in $\mathcal{T}$ share the edge. The abstract neighborhood $B(e)$ is the quotient of $\tilde{B}(e)$ obtained by identifying the preimages of the 1-singlex in each of the 3-simplex in $\triangle_{\mathcal{T}}$.

Let $\tilde{q}$ be a normal quadrilateral disk in a 3-simplex in $B(e)$. If it intersects the edge $e$, we assign a slope $s(\tilde{q}) \in \{\pm 1\}$ based on whether $\partial \tilde{q}$ crosses the equator of $B(e)$ from the Northern/Southern hemisphere to the Southern/Northern hemisphere, respectively. For each 1-singlex $e$ in $\mathcal{T}$, the *total slope* of a quadrilateral disk type $q \in \mathcal{T}$ with respect to $e$, written $s_e(q)$, is

$$s_e(q) := \sum_{\tilde{q} \text{ in } B(e) \text{ preimage of } q \text{ with } q \cap e \neq \emptyset} s(\tilde{q}).$$

If a quadrilateral disk type $q$ does not meet $e$, then define $s_e(q) = 0$. Note in a 3-simplex, there can be only two quadrilateral disk types meeting an edge, see Figure 14.

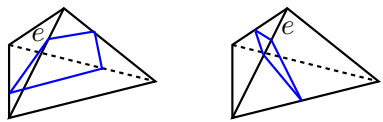

FIGURE 14. Quadrilaterals $q, q'$ meeting $e$.

Truncate the normal coordinates $\vec{F}$ (see (1) above) by removing all the coordinates for triangular disks to get $\vec{F}_Q$. Call $\vec{F}_Q$ the $Q$-coordinates of $\vec{F}$. For each singlex $e \in \mathcal{T}$, the $Q$-matching equations are the system of linear equations

$$\left\{ \sum_{i=1}^{t} (s_e(q_i') + s_e(q_i'') + s_e(q_i''')) = 0 \right\}_e.$$

Tollefson shows $\vec{F}_Q$ specifies a normal surface if and only if they satisfy the $Q$-matching equations at every edge $e$ of the triangulation $\mathcal{T}$.

**Theorem 2.3.** [Tol98, Theorem 1] *Let $M$ be a compact 3-manifold with a fixed triangulation $\mathcal{T}$. If $F$ is a normal surface in $M$ specified by $\vec{F}$ then the $Q$-coordinates $\vec{F}_Q$ give an admissible solution to the $Q$-matching equations. Moreover, if $\vec{v}$ is a nonzero admissible solution to the $Q$-matching equations, then there exists a unique normal surface $F$ in $M$ with no trivial components such that $\vec{F}_Q = \vec{v}$.*

2.5. **$Q$-matching equations for the diagrammatic triangulation $\mathcal{T}$.** In this section, we will write the $Q$-matching equations for the triangulation $\mathcal{T}$ constructed in the previous section in terms of the graph $G = G(D)$ as defined in Section 2.2. Fix such a graph $G$ corresponding to a knot diagram, and let $\mathcal{T}$ be the triangulation obtained in the previous section by applying the inflation procedure to the ideal triangulation $\mathcal{T}^*$ in Section 2.1. We say that a tetrahedron $\mathcal{T}$ is *associated* to a twist region $T$ if it comes from an octahedron placed at a crossing in the twist region, or if it is a tetrahedron in $\mathcal{T}$ added through the inflation procedure between faces of octahedra, each of which is placed at a crossing in the twist region. Denote the set of octahedra associated to a twist region by $O_T$. To describe the boundary faces of $O_T$, we will let it inherit the numbering of vertices from the corresponding faces of the first octahedron and the last octahedron of the twist region $T$. An edge of $\mathcal{T}$ is said to be *interior* to $O_T$ if it is not identified with any edges outside of the tetrahedra in $O_T$, or if it belongs to an edge shared by a tetrahedron in $O_T$ and a tetrahedron added through inflation.

**Definition 2.4.** We will denote by $R_G$ a complementary region, or, a face, of the graph $G$ in $S^2$. The unbounded complementary region containing $\infty$ in $S^2 = \mathbb{R}^2 \cup \infty$ is denoted by $R_\infty$.

Given a vertex $V$ of $G$ and an edge $E \in G$ which intersects $V$, let $O_T$ be the set of octahedra associated to the twist region $T$ represented by $E$ in $G$. We define the coordinates $C_i(E)$ for $i \in \{1, 2\}$ to be

$$C_1(E) := s_{O_T(34)}(q_i') + s_{O_T(34)}(q_i'') + s_{O_T(34)}(q_i''') \qquad C_2(E) := s_{O_T(12)}(q_i') + s_{O_T(12)}(q_i'') + s_{O_T(12)}(q_i''').$$

A twist region *intersects* $R_G$ if the corresponding edge in $G$ is a subset of $\partial R_G$. Given a face $R_G$ of $G$ and an edge $E \in G$ which intersects $R_G$, let $O_T$ be the set of octahedra associated to the twist region $T$ represented by $E$ in $G$. We define the coordinates $B_i(E)$ for $i \in \{1, 2\}$ to be, if $w(E) < 0$,

$$B_1(E) := s_{O_T(14)}(q_i') + s_{O_T(14)}(q_i'') + s_{O_T(14)}(q_i''') \qquad B_2(E) := s_{O_T(23)}(q_i') + s_{O_T(23)}(q_i'') + s_{O_T(23)}(q_i''').$$

**Lemma 2.5.** *Let $K$ be a nontrivial knot with diagram $D = D(K)$ and a planar connected weighted graph $G = G(D)$. Let $\mathcal{N} = \{\mathcal{N}_T\}$ be a collection of normal quadrilateral disks in the tetrahedra associated to each twist region $T$ of $G$ for the triangulation $\mathcal{T}$. Assume $\mathcal{N}$ satisfies the $Q$-matching equations at each interior edge in $\mathcal{T}$ of $O_T$. Then $\mathcal{N}$ represents a normal surface in $\mathcal{T}$ if*

*(a) for every vertex $V \in G$,* $\displaystyle \sum_{E \text{ intersecting } V} C_1(E) = \sum_{E \text{ intersecting } V} C_2(E) = 0$, *and*

*(b) for every face $R_G$ of $G$,* $\displaystyle \sum_{E \text{ intersecting } R_G} B_1(E) = \sum_{E \text{ intersecting } R_G} B_2(E) = 0$.

*Proof.* An edge of $\mathcal{T}$ is either an interior edge, where by assumption the $Q$-matching equations are satisfied, or an edge which is identified with an edge from another twist region around a vertex or around a complimentary region $R_G$. The statement of the lemma then follows from Theorem 2.3 applied to these edges. $\square$

2.6. **Normal surfaces local to a twist region.** For a twist region $T$ in a link diagram, let $c(T)$ denote the number of crossings in the twist region $T$. Fixing the twist region $T$, we number the crossings in $T$: $1, \ldots, j, \ldots, c(T)$. The octahedron associated to the $j$th crossing from $\mathcal{T}^*$ is denoted by $o_j$. The inflation procedure with our chosen frame $\Lambda$ adds a pair of tetrahedra between the $j$th and the $j+1$th crossing of the twist region $T$ which we will denote by $t_j$ and $t_j'$.

Assuming $T$ is a positive twist region, we consider three types of normal surfaces in the tetrahedra associated to $T$. They are determined by a choice of the three types of a quadrilateral disk $q_j', q_j'', q_j'''$ in the center tetrahedron $(o_j)(1234)$. See Figure 15.

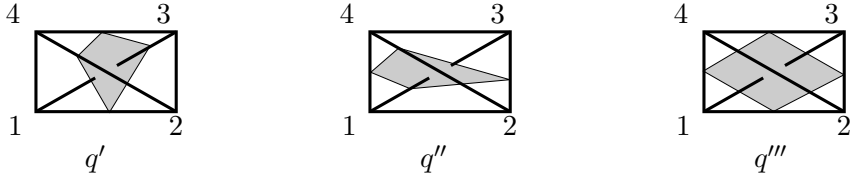

FIGURE 15. Three types of quads for the center tetrahedron.

**Local quad assignments**. For a negative twist region $T$ in a knot diagram fix $j \in \{1, \ldots, c(T)\}$. We specify the quad assignments for each tetrahedra associated to $T$. For now we assume that the transverse orientation along each branch is pointing into the tetrahedron added through inflation. For a positive twist region we take the mirror image of the quad assignments for that of the negative twist region. Define for a quad type $q$, $t(q)$ the index of the quad in the triangulation. Given an edge $e$ in $G$, a quad type $q$, and a tetrahedron $x$ in $\mathcal{T}$, let

$$n_{e,x}(q) := t(q') \text{ such that } q' \in x \text{ and } s_e(q) + s_e(q') = 0.$$

$$m_{e,x}(q) := t(q') \text{ such that } q' \in x \text{ and } s_e(q) + s_e(q') = \pm 2.$$

Let $k$ be an index for a tetrahedron in $\mathcal{T}$. Recall $t_j, t_j'$ denote the pair of tetrahedra added from inflation between $o_j$ and $o_{j+1}$.

TYPE I: Let

$$v_1', v_2', v_3', v_4' = n_{o_j(24), o_j(0124)}(q_j'), n_{o_j(13), o_j(1345)}(q_j'), n_{o_j(24), o_j(0234)}(q_j'), n_{o_j(13), o_j(1235)}(q_j'),$$

respectively. Define

$$\vec{v}_j'(k) := \begin{cases} 1 \text{ if } j \text{ is odd and } k \in \{t(q_j'), v_1', v_2'\} \\ 1 \text{ if } j \text{ is even and } k \in \{t(q_j'), v_3', v_4'\} \\ 0 \text{ otherwise.} \end{cases}$$

$$\vec{a}_j'(k) := \begin{cases} 1 \text{ if } k \in \{m_{o_j(15), t_j}(v_1'), m_{o_j(02), t_j'}(v_2')\} \\ 0 \text{ otherwise.} \end{cases}$$

TYPE II: Let

$$v_1'', v_2'', v_3'', v_4'' = n_{o_j(24), o_j(0234)}(q_j''), n_{o_j(13), o_j(1345)}(q_j''), n_{o_j(24), o_j(0124)}(q_j''), n_{o_j(13), o_j(1235)}(q_j''),$$

respectively. Define

$$\vec{v}_j''(k) := \left\{ \begin{array}{l} 1 \text{ if } j \text{ is odd and } k \in \{t(q_j''), v_1'', v_2''\} \\ 1 \text{ if } j \text{ is even and } k \in \{t(q_j''), v_3'', v_4''\} \\ 0 \text{ otherwise.} \end{array} \right.$$

$$\vec{a}_j''(k) := \left\{ \begin{array}{l} \vec{a}_j'(k) \\ 0 \text{ otherwise.} \end{array} \right.$$

TYPE III: Let

$$v_1''', v_3''' = n_{o_j(24), o_j(0124)}(q_j'''), n_{o_j(13), o_j(1235)}(q_j'''),$$

respectively. Define

$$\vec{v}_j'''(k) := \left\{ \begin{array}{l} 1 \text{ if } k \in \{t(q_j''')\} \\ 0 \text{ otherwise.} \end{array} \right.$$

$$\vec{a}_j'''(k) := \left\{ \begin{array}{l} 1 \text{ if } k \in \{v_1''', v_3'''\} \\ 0 \text{ otherwise.} \end{array} \right.$$

**$S_I^-(n, k, w)$** for $0 \le k < n$. We define $S_I^-(n, k, w)$ to be the normal surface represented by the vector

$$k \sum_{j=1}^{c(T)} \vec{v}_j' + (n-k) \sum_{j=1}^{c(T)} \vec{a}_j'.$$

**$S_{II}^-(n, k, r, w)$** for $0 \le r \le w$. Then $S_{II}^-(n, k, r, w)$ is defined to be the normal surface represented by the vector

$$k \sum_{j=1}^{c(T)} \vec{v}_{j,r}'' + (n-k) \sum_{j=1}^{c(T)} \vec{a}_{j,r}''.$$

**$S_{III}^-(n, k, w)$** for $0 \le r \le w$. Then $S_{III}^-(n, k, r, w)$ is defined to be the normal surface represented by the vector

$$k \sum_{j=1}^{c(T)} \vec{v}_{j,r}''' + (n-k) \sum_{j=1}^{c(T)} \vec{a}_{j,r}'''.$$

We define $S_I^+(n, k, w)$, $S_{II}^+(n, k, r, w)$, and $S_{III}^+(n, k, w)$ for a positive twist region $T_w$ with $w$ crossings by taking the mirror image of the assignments for $S_I^-(n, k, -w)$, $S_{II}^-(n, k, r, -w)$, and $S_{III}^-(n, k, -w)$.

**Lemma 2.6.** *Let $w < 0$. The assignments of quads $S_I^\pm(n, k, w)$, $S_{II}^\pm(n, k, r, w)$, $S_{III}^\pm(n, k, w)$ satisfy the Q-matching equation at each interior edge of the set of tetrahedra $O_T$ associated to a twist region $T$ with $-w$ number of crossings.*

*Proof.* We have three cases $S_I^-(n, k, w)$, $S_{II}^-(n, k, r, w)$, and $S_{III}^-(n, k, w)$ to consider, and the cases $S_I^+(n, k, w)$, $S_{II}^+(n, k, w)$, and $S_{III}^+(n, k, w)$ are analogous. By construction, because we specifically choose quads in tetrahedra which satisfy the $Q$-matching equations on the edges in which they meet, the conditions of Theorem 2.3 are satisfied at every interior edge of the set of tetrahedra associated to a twist region. $\square$

We conclude Section 2 with how to compute the boundary slopes of the normal surfaces constructed locally from our assignment of normal disks in each twist region.

**2.7. Slopes and boundary slopes of normal surfaces.** For a knot with diagram $D$, we compute the boundary slope of a normal surface $\mathcal{N}$ that restricts to a surface $\mathcal{N}_T$ defined in Section 2.6 in a twist region $T$. If $\mathcal{N}$ also satisfies the conditions of Lemma 2.5, then it is a normal surface. In this section, we compute the slope of such a normal surface $\mathcal{N}$ by summing over local contributions to the slope from each twist region.

Let $\Sigma$ be a smoothly embedded, connected, and properly embedded surface in $S^3 \setminus K$, no component of which is a 2-sphere or a surface parallel to a subsurface of $\partial(S^3 \setminus K)$. Suppose $\partial\Sigma \neq \emptyset$. Since $S^3 \setminus K$ is irreducible, all the components of $\partial\Sigma$ are mutually parallel non-trivial simple closed curves. Fix a framing of $K$ and let $\mu$ and $\lambda$ denote the meridian and longitude basis of $H_1(\partial(S^3 \setminus K); \mathbb{Z})$. All the components of $\partial\Sigma$ determine the same *slope* $p/q \in \mathbb{Q} \cup \infty$ such that $p[\mu] + q[\lambda]$ represents the homology class of $[\partial\Sigma] \in H_1(\partial(S^3 \setminus K); \mathbb{Z})$.

**Definition 2.7.** If the inclusion on fundamental groups

$$\iota_* : \pi_1(\Sigma) \hookrightarrow \pi_1(S^3 \setminus K)$$

induced by the inclusion $\iota$ of $\Sigma \subset S^3 \setminus K$ is injective, then the slope $p/q$ determined by the components of $\partial\Sigma$ is called a *boundary slope*.

Note that the slope $p/q$ is defined for any smoothly embedded, connected, and properly embedded surface in $S^3 \setminus K$, and we compute this quantity for $\mathcal{N}$. By construction, $\mathcal{N}$ is obtained by gluing the surfaces $\mathcal{N}_T$ in each twist region $T$. Each of these local surfaces are obtained by adding saddles to curves on the boundary. Thus, the boundary slope of the surface may be computed locally by summing over how the surface twists around the knot strands in a twist region. This uses the strategy employed in [HO89] for computing boundary slopes of Montesinos knots.

Passing through a saddle in the surface involves replacing one pair of opposite sides of a quadrilateral by the other pair of opposite sides. The change increases the slope of the surface if it twists through the clockwise direction. The change decreases the slope if it twists through the counter-clockwise direction. The total number $\tau(\mathcal{N})$ of twists is

$$\tau(\mathcal{N}) = 2(s_- - s_+)/m,$$

where $s_-/s_+$ is the number of slope decreasing/increasing saddles. The boundary slope is then

$$(2) \qquad\qquad s(\mathcal{N}) = \tau(\mathcal{N}) - \tau(\Sigma_0),$$

where $\Sigma_0$ is a Seifert surface from the Seifert algorithm applied to the diagram $D$ of the knot.

**Lemma 2.8.** *Suppose we have an assignment of normal disks $\mathcal{N} = \{\mathcal{N}_T\}$ for the triangulation $\mathcal{T}$ of $S^3 \setminus K$ that satisfies conditions (a) and (b) of Lemma 2.5, and each $\mathcal{N}_T$ is the surface $S_I^{\pm}(n, k, w)$, $S_{II}^{\pm}(n, k, w)$, or $S_{III}^{\pm}(n, k, w)$ defined in Section 2.6. The contributions of a Type I, Type II, or Type III surface to the twist number of $\mathcal{N}$ are tabulated as follows:*

| Surface | $S_I^-(n,k,w)$ | $S_{II}^-(n,k,r,w)$ | $S_{III}^-(n,k,w)$ |
|---------|----------------|---------------------|--------------------|
| $\tau$ | $-2(1-\frac{k}{n})$ | $2(w-r+1-\frac{k}{n})$ | $0$ |
| Surface | $S_I^+(n,k,w)$ | $S_{II}^+(n,k,w)$ | $S_{III}^+(n,k,w)$ |
| $\tau$ | $2(1-\frac{k}{n})$ | $-2(w-r+1-\frac{k}{n})$ | $0$ |

TABLE 4. The local contribution to the slope of the surface represented by the normal subset $\mathcal{N}$.

*Proof.* We count the number of twists that the surfaces make along the boundary of a twist region. $\qquad\square$

**Lemma 2.9.** *The twist number of a normal subset $\mathcal{N} = \{\mathcal{N}_T\}$ that satisfies conditions* (a) *and* (b) *of Lemma 2.5, and where each $\mathcal{N}_T$ is the surface $S_I^{\pm}(n, k, w)$, $S_{II}^{\pm}(n, k, w)$, or $S_{III}^{\pm}(n, k, w)$ defined in Section 2.6, is given by summing over all the local contributions:*

$$\tau(\mathcal{N}) = \sum_{T\,a\,twist\,region} \tau(\mathcal{N}_T).$$

*Proof.* We do not introduce changes in the homology class of the boundary curves when the surfaces in each set of tetrahedra associated to a twist region are identified with each other. Thus, the overall twist number of $\mathcal{N}$ is the sum over the twist number of the individual surface in each twist region.                                                                                      □

## 3. Colored Khovanov homology

In this section we describe colored Khovanov homology, which categorifies the colored Jones polynomial following the conventions of [BN07], [CK12]. We will define colored Kauffman states relate them to the normal surfaces defined in the previous section to prove Theorem 1.1.

### 3.1. The Temperley-Lieb algebra and Jones-Wenzl projectors.
For a fixed $n$, the $n$th Temperley-Lieb algebra $TL_n$ is a formal vector space of oriented link and tangle diagrams in the 2-disk $(\mathcal{D}^2, n)$ with coefficients in $\mathbb{Z}[q, q^{-1}]$, visualized as a rectangle with $n$ marked points on the top and bottom boundary, and modulo the Kauffman bracket skein relations.

- $\langle \diagbox \rangle = q \langle \rangle \langle \rangle - q^2 \langle \smallsmile \rangle$         $\langle \diagbox \rangle = q^{-2} \langle \smallsmile \rangle - q^{-1} \langle \rangle \langle \rangle$.
- $\langle \bigcirc \rangle = (q + q^{-1}) \langle \rangle$.

Elements of $TL_n$ are called *skein elements*. The multiplication operation $\cdot$ making $TL_n$ into an algebra sends two skein elements $\mathcal{U}, \mathcal{V} \in TL_n$ to the skein element $\mathcal{U} \cdot \mathcal{V}$ obtained by stacking the square containing $\mathcal{U}$ on top of the square containing $V$, and identifying the $n$ boundary points. The identity of the algebra operation, denoted by $1_n$, is $n$ parallel strands connecting the $n$ boundary points. A presentation of $TL_n$ is given by standard generators $1_n$ and $\{e_i\}, 0 < i < n$.

The *Jones-Wenzl projectors* $p_n \in TL_n$ are idempotent elements in the algebra that are uniquely defined by the following properties.

(1) $p_n - 1_n$ belongs to the subalgebra generated by $\{e_1, e_2, \dots, e_{n-1}\}$,
(2) $p_n \cdot p_n = p_n$, and
(3) $e_i \cdot p_n = p_n \cdot e_i = 0$ for all $1 \leq i \leq n - 1$.

The *quantum integer* $[n]$ is defined to be $[n] = \frac{q^n - q^{-n}}{q - q^{-1}}$.

The projectors can be defined, as when they first appeared in [Wen87], by the recurrence relation

$$p_1 = 1_n, \qquad p_n = p_{n-1} \sqcup 1_n - \frac{[n-1]}{[n]} p_{n-1} \cdot e_{n-1} \cdot p_{n-1}.$$

We will depict $p_n$ as a box on $n$ strands. Pictorially, the recurrence relation is then as shown in Figure 16.

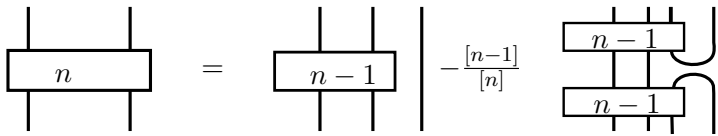

FIGURE 16. Recurrence relation for the Jones-Wenzl projector.

## 3.2. Bar-Natan's formulation [BN07] of Khovanov homology for tangles and cobordisms.
Letting the disk $\mathcal{D}^2$ be identified with the square $[0,1] \times [0,1]$ with $n$ marked points on $\{0,1\} \times [0,1]$, let $\mathrm{Cob}_n$ be the category with

- **Objects:** Isotopy classes of $q$-graded 1-submanifolds $T$ (arcs and circles) in $\mathcal{D}^2$ that are properly embedded (ie, $\partial T \subset \partial \mathcal{D}^2$), such that the endpoints of the arcs form a subset of the $n + n = 2n$ marked points on the boundary of $\mathcal{D}^2$.
- **Morphisms:** Isotopy classes of embedded cobordisms in $\mathcal{D}^2 \times [0,1]$ between objects of $\mathrm{Cob}_n$, which may or may not be decorated with dots, considered up to isotopy fixing the boundary, and modulo the following local relations.

(1) 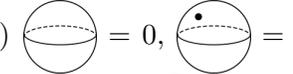 $= 0,$ 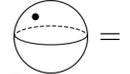 $= 1,$ 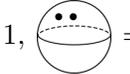 $= 0.$ In particular, let $S$ be an embedded cobordism, then $S \sqcup$ 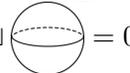 $= 0,$ $S \sqcup$ 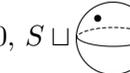 $= S,$ and $S \sqcup$ 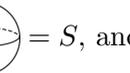 $= 0.$

(2) 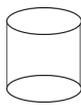 $=$ 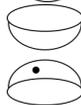 $+$ 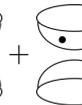 , which implies

(3) 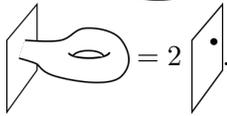 $= 2$ 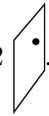 .

We refer to these relations as the *Bar-Natan* skein relations. The degree of a cobordism $C : q^i A \to q^j B$ is given by

$$\deg(C) = \deg_t(C) + deg_q(C),$$

where the topological degree $\deg_t(C) = \chi(C) - n$ is given by the Euler characteristic of $C$ and the $q$-degree $\deg_q(C) = j - i$ is given by the relative difference in the $q$-gradings. The maps $C$ used throughout the paper will satisfy $\deg(C) = 0$.

We make $\mathrm{Cob}_n$ into an additive category by formally adding finite direct sums of objects, and denote the result by $\mathbf{Cob}_n$.

The skein relation becomes

(3) 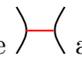 $= q$ $\Big)$ $\Big($ $\longrightarrow$ $q^2$ 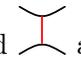 $\quad$ $= q^{-2}$ $\longrightarrow$ $q^{-1}$ $\Big)$ $\Big($ ,

where $\rangle\!\!-\!\!\langle$ and $\asymp$ are saddle cobordisms.

**Definition 3.1.** Let $\mathbf{Kom}_n$ be the category of chain complexes of finite direct sums of objects in $\mathbf{Cob}_n$. We allow chain complexes $C_*$ of unbounded positive homological degree, and require that for each chain complex $C_*$ there exists $N \in \mathbb{Z}$ such that $C_n = 0$ for all $n < N$.

To a tangle or a link diagram, the skein relation (3) associates an object in $\mathbf{Kom}_n$. Thus, an object in a chain complex also has a homological grading which we denote by $h$, so $\deg_h(B) = h(B)$. If $A \to B$ then $h(B) = h(A) + 1$. The category $\mathbf{Kom}_n$ is a categorification of the Templerley-Lieb algebra $TL_n$ in the sense that the Grothendieck group $K_0(\mathbf{Kom}_n)$ of $\mathbf{Kom}_n$ is isomorphic to $TL_n$ as a $\mathbb{Z}[q^{-1}][[q]]$-algebra [CK12, Lemma 2.10].

For a link diagram, where applying the skein relation at each crossing results in a set of disjoint closed curves, one then obtains Khovanov homology by applying the functor to $\mathbf{Kom}_0$ that sends the category $\mathbf{Cob}_0$ to to $\mathbb{Z}Mod$, the category of graded $\mathbb{Z}$-modules which maps disjoint unions to tensor products.

**Definition 3.2.** Let $V$ be the graded $\mathbb{Z}$-module freely generated by two elements $v_\pm$ with $\deg_q(v_\pm) = \pm 1$. Let $\mathcal{F}$ be the functor defined by

$$\mathcal{F}(\ \bigcirc\ ) = V$$
$$\mathcal{F}(\ \smallsmile\ ) = \epsilon : \mathbb{Z} \to V; \quad \epsilon : 1 \mapsto v_+$$
$$\mathcal{F}(\ \frown\ ) = \eta : V \to \mathbb{Z}; \quad \eta : v_+ \mapsto 0, \quad v_- \mapsto 1$$
$$\mathcal{F}(\ \asymp\ ) = s : V \to V \otimes V; \quad v_+ \mapsto v_+ \otimes v_- + v_- \otimes v_+, \quad v_- \mapsto v_- \otimes v_-$$
$$\mathcal{F}(\ \bigcirc\!\!-\!\!\bigcirc\ ) = m : V \otimes V \to V; \quad v_+ \otimes v_- \mapsto v_-, \quad v_- \otimes v_+ \mapsto v_-, \quad v_+ \otimes v_+ \mapsto v_+,$$
$$v_- \otimes v_- \mapsto 0.$$

3.3. **Categorification of the Jones-Wenzl projectors.** Several constructions have been given for a categorification of the Jones-Wenzl projector, see [CK12], [Roz14a], [Hog19], [Hog20], and [Kho05]. To summarize, the goal is to construct a chain complex which decategorifies to the Jones-Wenzl projector in $TL_n$. Let $K_0(\mathbf{Kom}_n)$ be the Grothendieck group of $\mathbf{Kom}_n$. For the construction in [CK12], this means that an object $P_n$ of $\mathbf{Kom}_n$ is associated to the projector in $TL_n$, such that the image of $P_n$ in the isomorphism between $K_0(\mathbf{Kom}_n)$ and $TL_n$ is the projector $p_n$. The categorification in [CK12], [Roz14a], and [Hog19], [Hog20] are unique up to homotopy since they satisfy the categorified versions of the identities characterizing the Jones-Wenzl projector described below, except for the construction by Khovanov in [Kho05][1]:

**Definition 3.3** ([CK12, Definition 3.1]). A chain complex $(P_*, d_*) \in \mathbf{Kom}_n$ is a *universal projector* if

(1) It is positively graded with degree zero differential.
   (a) $P_k = 0$ for all $k < 0$ and $\deg_q(P_k) \geq 0$ for all $k > 0$.
   (b) $d_k$ is a matrix of degree zero maps for all $k \in \mathbb{Z}$.
(2) The identity diagram appears only in homological degree zero and only once
   (a) $P_0 \cong 1$.
   (b) $P_k \not\cong 1 \bigoplus D$ for any $D \in Mat(Cob(n))$ for all $k > 0$.
(3) The chain complex $P_*$ is contractible when composed with turnbacks. That is for any generator $e_i \in TL_n$, $0 < i < n$,
   (a) $P_* \otimes e_i \simeq 0$
   (b) $e_i \otimes P_* \simeq 0$

If a universal projector satisfying Definition 3.3 exists, then it is unique up to homotopy, and satisfies the categorical analogues of the defining properties of the Jones-Wenzl projector in the Temperley-Lieb algebra described in Section 3.1.

**Theorem 3.4** ([CK12, Corollary 3.5]). *If $\mathcal{C}, \mathcal{C}' \in Kom(n)$ are universal projectors then $\mathcal{C} \simeq \mathcal{C}'$.*

From this, one sees that the constructions of [CK12], [Roz14a], [Hog19], and [Hog20] are equivalent.

With the product $\cdot$ in $TL_n$ replaced by the tensoring operation $\otimes$ of complexes, a universal projector $P$ satisfies the following properties:

(i) $P_* \otimes P_* \simeq P_*$.
(ii) $P_* \otimes e_i \simeq 0 \simeq e_i \otimes P_*$.

For the categorification of the colored Jones polynomial we will be using the version in [CK12], denoted by $P_n$, which constructs the $n$th Jones-Wenzl projector from the $n-1$th projector by categorifying the following identity:

---

[1]We would like to remark that Khovanov's categorification of the colored Jones polynomial [Kho05] is known not to be isomorphic to the constructions by [CK12].

(4)

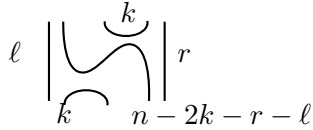

FIGURE 17. Recursion relation for the Jones-Wenzl projectors.

For each $n$ we will denote the $n$ terms on the right hand side of the above recurrence relation expanding the $n$ projector as $p_{n,i}$, where $p_{n,1} = p_{n-1} \sqcup 1$, and $p_{n,i}$ is the element with bottom cap joining the points $n - i + 1$ to $n - i + 2$ for $i > 1$. The corresponding coefficient functions in $q$ are denoted by $f_{n,i}$.

The complex $P_n$ is constructed recursively from $P_{n-1}$ as depicted below, and becomes periodic with period $2(n-1)$ after the first $2n$ terms.

$$p_{n,1} \to qp_{n,2} \to \cdots \to q^{n-1}p_{n,n} \to q^{n+1}p_{n,n} \to q^{n+2}p_{n,n-1} \to \cdots \to q^{2n-1}p_{n,2} \to$$
$$\to q^{2n+1}p_{n,2} \to \cdots$$

The categorification complex $P_n$ consists of crossingless skein elements in $TL_n$, where the identity skein element $|_n$ in $TL_n$ has homological degree 0. We will denote by $\overline{\deg_h}(\mathcal{J})$ the smallest possible homological degree of a skein element.

**Lemma 3.5.** *Assume $\ell, r, k, n$ are nonnegative integers with $\ell + 2k + r \leq n$. The skein element $\mathcal{J}_{\ell,r,k}^n$ in Figure 18 in the colored Khovanov complex has homological degree $\overline{\deg_h}(\mathcal{J}_{\ell,r,k}^n) = k(n - r - k - \ell)$.*

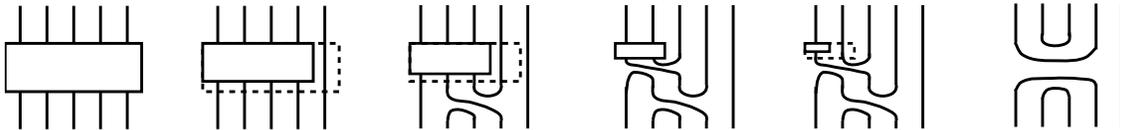

FIGURE 18. The skein element $\mathcal{J}_{\ell,r,k}^n$.

*Proof.* We describe a sequence of choices of terms, starting with $p_n$ in the recurrence relation of Figure 17 to arrive at $\mathcal{J}_{\ell,r,k}^n$. First choose $p_{n,1}$, then for the Jones-Wenzl projector $p_{n-1}$ in $p_{n,1}$ choose the diagram $p_{n-1,1}$. Repeat with the Jones-Wenzl projector $p_{n-2}$, until we choose $p_{n-r+1,1}$ to arrive at the skein element $p_{n-r} \sqcup 1_r$, then choose $p_{n-r,n-r-k-\ell+1}$ to replace the Jones-Wenzl projector $p_{n-r}$. This is followed by the choices of $p_{n-r-1,n-r-k-\ell+1}, \ldots, p_{n-r-k+1,n-r-k-\ell+1}$. For the last part of the sequence we just choose $p_{n-r-k}, 1, \ldots p_2, 1$ until there is no longer a projector. See Figure 19 for an example.

FIGURE 19. A sequence of choices of expansions for the projector.

Each of these choices changes the homological degree as follows: For $p_{n,i}$, the homological degree increases by $i - 1$. Summing this over the sequence described above, where in particular, the terms $p_{*,1}$ do not contribute to any change in the homological grading, we get the homological degree $= k(n - r - k - \ell)$. $\square$

3.4. **A definition of colored Khovanov homology.** Let $D$ be an oriented link diagram, the colored Khovanov chain complex $C\mathcal{KH}^n(D)$ categorifying the $n$ colored Jones polynomial $J_K^n(q)$ is obtained by tensoring the chain complex $P_n$ categorifying $p_n$ with the $n$-blackboard cable of the link diagram[2]. Colored Khovanov homology is then the homology groups $\mathcal{KH}^n$ of the chain complex $C\mathcal{KH}^n(D)$. See Figure 20 for an illustration.

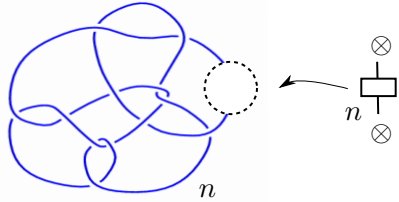

FIGURE 20.  Forming the complex for $C\mathcal{KH}^n$ by tensoring.

The colored Jones polynomial can be recovered from colored Khovanov homology through decategorificaton. That is, taking the graded Euler characteristic of the homology groups $\mathcal{KH}^n$.

$$J_K^n(q) = \sum_{i,j} (-1)^i q^{i+j} \dim(\mathcal{KH}_{i,j}^n(K)).$$

Here $i$ is the homological grading and $j$ is the quantum grading. With our convention, the colored Jones polynomial of the unknot is

$$J_K^n(q) = (q + q^{-1})^n.$$

3.5. **Colored Kauffman states.** We will study the colored Khovanov complex from the skein element $D^n$, where $D^n$ is obtained from the knot diagram $D$ by $n$-cabling a strand with a Jones-Wenzel projector. Using Theorem 3.4 (i), we double the projector, so that every $n$-cabled twist region corresponding to an edge $E$ in the graph $G(D)$ (from Section 2.2) is framed by four projectors.

The next two results allow us to decompose the chain complex as a direct sum of crossingless colored tangles decorated by Jones-Wenzl projectors.

**Theorem 3.6.** [Roz14b, Theorem 3.5] *The colored Khovanov chain complex of the crossing of two $n$ colored strands can be presented as a multi-cone of crossingless colored tangles $\mathcal{J}(k)$, see Figure 21, for $0 \leq k \leq n$, with*

$$\overline{\deg_h}(\mathcal{J}(k)) = n^2 - k^2 \qquad \text{if the crossing is positive,}$$

*and*

$$\overline{\deg_h}(\mathcal{J}(k)) = k^2 \qquad \text{if the crossing is negative.}$$

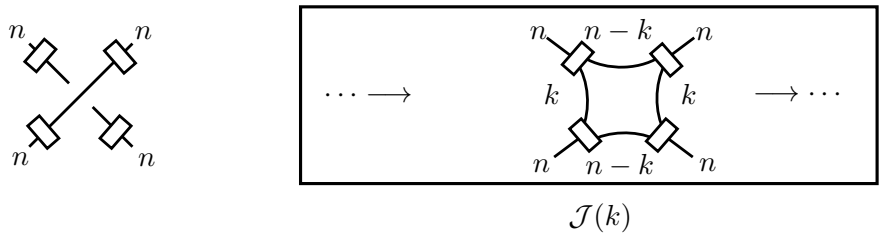

FIGURE 21.  A multi-cone presentation at a colored crossing.

---

[2]The categorification would be different if it is constructed via the methods of [Kho05]. In particular the complexes he constructs are bounded in terms of homological width, whereas the construction in [CK12] is unbounded in positive homological degree.

By expanding the Jones-Wenzl projectors appearing in the complex of an $n$-cabled twist region framed by four Jones-Wenzl projectors after applying Theorem 3.6 via (4), we also get an expansion of the colored chain complex as a multi-cone of crossingless colored tangles $\mathcal{S}(k)$ over twist regions.

**Corollary 3.7.** *The colored Khovanov chain complex of a twist region with $w$ crossings of two $n$ colored strands can be presented as a multi-cone of crossingless colored tangles $\mathcal{S}(k)$, see Figure 22, for $0 \le k \le n$, with*

$$\overline{\deg_h}(\mathcal{S}(k)) = w(n^2 - k^2) \qquad \text{if the twist region is positive,}$$

*and*

$$\overline{\deg_h}(\mathcal{S}(k)) = wk^2 \qquad \text{if the twist region is negative.}$$

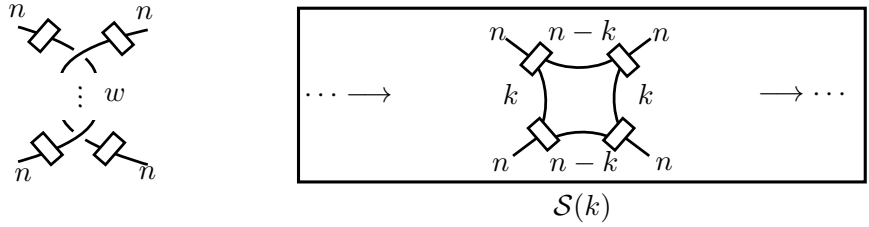

FIGURE 22. A multi-cone presentation at a colored twist region.

**Definition 3.8.** Number the twist regions of $D$ corresponding to edges of the graph $G = G(D)$ $1, \ldots, |E(G)|$, and fix $n$. A *colored Kauffman state* $(\sigma, \mathcal{J})$ of $\mathcal{CK}\langle^n(D)$ is a map $\sigma$ on the set of twist regions $\Gamma(D)$ of $D$.

$$\sigma : \Gamma(D) \overset{Id}{\to} E(G) \to \{0, 1, \ldots, n\},$$

along with a map $\mathcal{J}$ from the set of Jones-Wenzl projectors into the set of crossingless matchings in the expansion of the $n$th Jones-Wenzl projector via the recursion relation (4).

The *parameter vector* $k(\sigma)$ of a colored Kauffman state $\sigma$ is the vector recording the value of $\sigma$:

$$k(\sigma) := (\sigma(1), \ldots, \sigma(|E(G)|)).$$

Applying a colored Kauffman state to a link diagram $D$ means replacing an $n$-cabled twist region in $D^n$ with the skein element $\mathcal{J}(k)$, and the Jones-Wenzl projectors replaced by an term in its expansion. This will result in a set of disjoint closed curves as when applying a Kauffman state to a diagram.

### 3.6. **Colored surface states from flows on graphs.**

**Definition 3.9.** [MOY98][3] Let $G$ be an oriented planar graph. A *flow $f$ on $G$* is a map from the edge set of $G$ to positive integers less than or equal to $n$ such that for every vertex $v \in G$ the sum of its values on the edges coming into $v$ is equal to that on the edges going out from $v$.

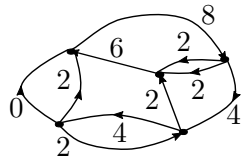

FIGURE 23. A flow on a directed graph $G$ indicated by the blue numbers.

---

[3]Note [MOY98] requires the graph to be trivalent but we do not require that here.

Let $(\sigma, \mathcal{J})$ be a colored Kauffman state on $D^n$. Applying the state to $D^n$ means that we replace the $n$-cable of each twist region $T$ with $\sigma(T)$, and we replace each Jones-Wenzl projector by its image under $\mathcal{J}$. We will call the resulting set of disjoint closed curves in the plane the *state circles* of the colored Kauffman state.

Now we let $G(D)$ be a weighted planar graph coming from a knot diagram $D$ as in Section 2.2. Fix an edge $E \in G$. A state circle which intersects $\sigma(T_E)$ in an arc determines a finite walk from $v_o$ to $v_t$ with the other arc making up the circle. We will consider colored Kauffman states for whom every such state circle for every edge in $G$ actually determines a finite *path* from $v_0$ to $v_t$.

In addition, we orient the edges of the graph $G$ based on an orientation of these paths from $v_o$ to $v_t$ from the state circles of the colored Kauffman state, such that two paths through the same edge share the same orientation. If an edge $E$ of $G$ intersecting a path is such that $w(E) > 0$, then the edge $E$ has the same orientation as that of the path. Otherwise, it has the opposite orientation.

**Definition 3.10.** Given a colored Kauffman state $(\sigma, \mathcal{J})$ on $D^n$, label the edge $E$ in $G$ with $\sigma(T_E)$ and consider the resulting map $f : E(G) \to \{1, \ldots, n\}$. We say that $\sigma$ *induces* a flow on $G$ if $f$ is a flow.

**Definition 3.11.** Given a colored Kauffman state $(\sigma, \mathcal{J})$ inducing a flow, we say that $(\sigma, \mathcal{J})$ is *compatible* if the expansions chosen by $\mathcal{J}$ preserve the flow induced by $\sigma$ on $G$.

**Lemma 3.12.** *For every colored Kauffman state $\sigma$ inducing a flow $f$, there is a choice of expansions of Jones-Wenzl projectors $\mathcal{J}$, such that $(\sigma, \mathcal{J})$ is compatible.*

*Proof.* Fix a vertex $v$ in the graph $G$. Condition 3.13 (2) forces there to be only one inflow or outflow. By Lemma 3.5, we may choose the generator in the expansion of the Jones-Wenzl projector giving the requisite number of through strands. □

Consider the graph $G_-$ obtained by removing all the positive edges. A negative edge $E$ is said to be *disjoint* from another negative edge $E'$ if $E$ and $E'$ belong to different components in $G_-$.

To reduce the technicality of the proof, we will consider colored Kauffman states $\sigma$ inducing nonzero flows that satisfy the following additional conditions:

**Condition 3.13.**
   (1) At every vertex of the graph $G$, $\displaystyle\sum_{E \text{ intersects } v, w(E)>0} k_E \leq n$.
   (2) If there are any edges $E$ with $w(E) < 0$ which intersects $v$, then all but one has $k_E = 0$.
   (3) At each vertex $v$, either there is a single edge with nonzero flow into $v$, or there is a single edge with nonzero flow out of $v$.
   (4) If there is a negative edge $E$ intersecting $v$, then $k_E = n$.
   (5) A state circle of the colored Kauffman state does not connect distinct negative components. It does not connect a negative twist region to itself, nor contain an edge corresponding to a negative twist region in a different component.

With the orientation on $G$ from a colored Kauffman state $\sigma$ inducing a flow, we assign normal disks constructed in Section 2.6 to a colored Kauffman state on $D^n$:
If the flow through the twist region $T$ belongs to the same path in $G$, then

$$\mathcal{N}_T(\sigma) = \begin{cases} S_I^{\pm}(n, k, w) \text{ if } w(E) = \pm 1 \text{ and } 0 \leq \sigma(T) < k_E \\ S_{II}^{\pm}(n, k, r, w) \text{ if } w(E) = \pm 1 \text{ and } \sigma(T) = n, \end{cases}$$

Otherwise,

$$\mathcal{N}_T(\sigma) = S_{III}^{\pm}(n, k, w) \text{ if } k_E = \pm 1.$$

We show when this assignment of normal disks to a colored Kauffman state $\sigma$ gives a normal surface $\mathcal{N}_\sigma = \{\mathcal{N}_T(\sigma)\}$ by translating the conditions of Lemma 2.5. Let $R_G$ be a face of $G$, define

$$B_{R_G}(E, \sigma) = \left\{ \begin{array}{l} s_{O_T(14)}(\mathcal{N}_T) \text{ if } O_T(14) \text{ intersects } R_G \\ s_{O_T(23)}(\mathcal{N}_T) \text{ if } O_T(23) \text{ intersects } R_G \end{array} \right.$$

**Proposition 3.14.** *Let $\sigma$ be a colored Kauffman state that induces a flow on $G$ satisfying Condition 3.13. Then if the parameters $k = k(\sigma)$ satisfies the following criteria:*

(1) *Around a face $R_G$ of $G$,* $\displaystyle\sum_{E \text{ intersects } R_G} B_{R_G}(E, \sigma) = 0$

(2) *At every vertex of the graph $G$,* $\displaystyle\sum_{E \text{ intersects } v} sgn(E)k_E = 0$. *Here $sgn(E)$ is the sign of $w(E)$.*

*Then $\mathcal{N}_\sigma$ defines a normal surface in the triangulation $\mathcal{T}$.*

*Proof.* This follows from translating conditions (a) and (b) of Lemma 2.5 for the assignment of the normal surface $\mathcal{N}_\sigma$ to colored Kauffman state. $\square$

We will first enlarge the assignment of surfaces to a colored Kauffman state to state surfaces.

A Kauffman state $\sigma$ on a link diagram $D$ gives a surface $S_\sigma$ as follows. Each state circle of $\sigma$ bounds a disk in $S^3$. This collection of disks can be disjointly embedded in the ball below the projection plane. At each crossing of $D$, we connect the pair of disks neighboring the crossing by a half-twisted band to construct a surface $S_\sigma \subset S^3$ whose boundary is $D$.

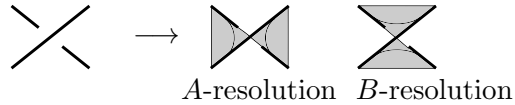

$A$-resolution $\quad$ $B$-resolution

FIGURE 24. Twisted bands at a crossing.

**Definition 3.15.** A surface $S_\sigma$ thus obtained from a Kauffman state $\sigma$ on a link diagram is called a *state surface*.

**Definition 3.16.** A Kauffman state $\sigma$ is *$\sigma$-adequate* if every edge in the state graph has two ends on distinct state circles.

Suppose $\widetilde{\sigma}$ is a Kauffman state on a link diagram $D$. Let $(\sigma, \mathcal{J} \equiv Id)$ be a colored Kauffman state which chooses the same resolution on the $n^2$ cabled crossings in $D^n$ corresponding to a crossing in $D$ as $\widetilde{\sigma}$ on $D$, and chooses the identity expansion for every Jones-Wenzl projector.

**Definition 3.17.** Fix $n$. A colored Kauffman state $\sigma$ is a *colored surface state* if it satisfies the conditions of Proposition 3.14 and Conditions 3.13, or if it comes from a $\widetilde{\sigma}$-adequate Kauffman state $\widetilde{\sigma}$ on $D$.

The reason for restricting to $\sigma$-adequate Kauffman states on $D$ is that a non-$\sigma$ adequate state is not essential since an edge in the state graph has two ends on the same state circle.

If a colored surface state $\sigma$ satisfies the conditions of Proposition 3.14, then we assign the surface $\mathcal{N}_\sigma$ to $\sigma$. If the colored surface state $\sigma$ comes from a Kauffman state $\widetilde{\sigma}$ on $D$, then we assign the surface $S_{\widetilde{\sigma}}$. Note this is well-defined since if a colored surface state satisfies the conditions of Proposition 3.14 and also comes from a Kauffman state $\sigma$, then the assigned surfaces coincide.

3.7. **Recovering the slope of a normal surface.** We will denote by $\mathcal{N}_0$ the surface obtained from the colored surface state $\sigma_0$ which has $k(\sigma_0) = \vec{0}$, and call it the *reference surface*. The reference surface exists for every nontrivial knot and is always a spanning surface. However, it is not always essential [Oza11].

Let $\sigma$ be a colored surface state and fix $n$. Write $h(\sigma) = a_\sigma n^2$ and note $a_{\sigma_0} = s(\mathcal{N}_0)$ by [FKP13].

**Lemma 3.18.** *If $a_\sigma - a_{\sigma_0} = \tau(\mathcal{N}) - \tau(\mathcal{N}_0)$, then $a_\sigma = s(\mathcal{N})$.*

*Proof.* Let $S$ be the Seifert surface obtained by applying the Seifert algorithm to the diagram $D$. We know

$$s(\mathcal{N}) = \tau(\mathcal{N}) - \tau(S).$$

By assumption,

$$s(\mathcal{N}) = a_\sigma - a_{\sigma_0} + \tau(\mathcal{N}_0) - \tau(S)$$

Then

$$s(\mathcal{N}) = a_\sigma - (\tau(\mathcal{N}_0) - \tau(S)) + \tau(\mathcal{N}_0) - \tau(S) = a_\sigma.$$

$\square$

**Theorem 1.1.** *Let $K$ be a nontrivial knot in $S^3$ with diagram $D$. Suppose a colored Kauffman state $\sigma$ is a colored surface state with homological grading $h_\sigma$, then there is a corresponding normal surface $\mathcal{N}_\sigma$ in the octahedral triangulation $\mathcal{T}_D$ of the knot $K$ with slope $s_\sigma$, such that*

$$h_\sigma = s_\sigma n^2.$$

*Proof.* That $\mathcal{N}_\sigma$ is a normal surface follows directly from Proposition 3.14. The twist number of the normal surface is the sum over the local contribution of each surface to the twist number. This is detailed in Lemma 2.8.

Suppose the colored surface state has parameter $k(\sigma) = \vec{0}$. Note this is always a solution and gives a state surface. One computes the boundary slope of the state surface using [FKP11] and this gives the statement in the theorem for this case.

Otherwise suppose $k(\sigma) \neq 0$. The state circles in the diagram after applying the colored Kauffman state can be partitioned into unique paths $p_1, \ldots, p_s$ with corresponding edges $E_1, \ldots, E_2$ in $G$. The homological degree $h_\sigma$ is written, summing over all the twist regions using Corollary 3.7, as

$$h_\sigma = \sum_{i=1}^s \left( \left( \sum_{E \in p_i} w(E) \right) - 1 \right) k_{E_i}^2.$$

On the other hand, Condition 3.13 guarantees that if $E_i$ and $E_j$ share a complementary region $R_G$, then

$$\left( \left( \sum_{E \in p_i} w(E) \right) - 1 \right) k_{E_i} = \left( \left( \sum_{E \in p_j} w(E) \right) - 1 \right) k_{E_j}.$$

Thus if we further partition $p_1, \ldots, p_s$ according to the common connected negative component that they share, then

$$h_\sigma = 2 \left( \left( \sum_{E \text{ in negative component}} (-w_E + 1) \right) + (q_1 - 1)\frac{k_{E_1}}{n} \right) n^2,$$

and $\sum_i k_{E_i} = n$. Finally,

$$a_\sigma - a_{\sigma_0} = 2 \left( \left( \sum_{E \text{ in negative component}} (-w_E + 1) \right) + (q_1 - 1)\frac{k_{E_1}}{n} \right) = \tau(\mathcal{N}) - \tau(\mathcal{N}_0)$$

by Lemma 2.8. The theorem follows from Lemma 3.18.                                    $\square$

3.8. **Homology classes of colored state surfaces.** In this section we specify precisely the generator in colored Khovanov homology from a colored surface state whose homological grading corresponds to the slope of the associated surface in the sense of Theorem 1.1 and we prove Theorem 1.2.

We define the generator corresponding to the colored Kauffman state with a compatible expansion of the Jones-Wenzl projector as $X_\sigma$ corresponding to the element $v_- \otimes \cdots \otimes v_-$.

**Theorem 1.2.** *Suppose $K$ is a nontrivial knot with a highly twisted diagram $D$. For a colored surface state $\sigma$ with $n > 1$ suppose that the associated surface $\mathcal{N}_\sigma$ is essential. Then the corresponding generator $X_\sigma$ in colored Khovanov homology is a cycle and $[X_\sigma] \neq 0$.*

*Proof.*
**Case 1:** $S_\sigma$ is a state surface. In this case the corresponding Kauffman state must be $\widetilde{\sigma}$-adequate. Otherwise $S_\sigma$ would be essential because of the existence of a one-edged loop. Similarly, $\sigma$ is also $\sigma$-adequate. This means that $d(X_\sigma)$ merges a pair of circles corresponding to $v_- \otimes v_-$, and so $d(X_\sigma) = 0$. This proves $X_\sigma$ is a cycle.

When $n = 1$. The homology class $[X_\sigma]$ is nonzero since every circle in $X_\sigma$ that results from splitting in $d^{-1}$ is marked $+$. We are done if there is one way to split into the circle (therefore, a single saddle from state with homological grading $h - 1$ that splits the given circle). If there are multiple ways to split this particular circle, then we must have a set of parallel edges on a pair of state circles. Canceling out in pairs gives us the result. When $n > 0$, all the previous states would have a loop composed with a cap/cup which would go to zero.

**Case 2:** $S_\sigma$ is not a state surface. This means that it has multiple sheets and comes from a case where the vector $k$ is not zero and all the circles are marked with a $-$, meaning that the generator is equal to $v_- \otimes \cdots \otimes v_-$. We show that the corresponding generator $X_\sigma$ gives a homology class by computing its boundary map. Because the diagram is highly twisted, changing the resolution at a single crossing or changing the expansion at a single Jones-Wenzl projector merges a pair of circles. Since all the circles are marked with a $-$, its image under the boundary map goes to 0. This shows that $X_\sigma$ gives a cycle in colored Khovanov homology.

Next we show $[X_\sigma] \neq 0$. Fix a negative twist region, we consider the decomposition of the chain complex **CKh** where we expand all the Jones-Wenzl projector except for those four framing the negative twist region, $\widetilde{X_\sigma}$. The element $X_\sigma$ is contained in the expansion of the Jones-Wenzl projectors $\widetilde{X_\sigma}$ where all four of the remaining projectors are replaced by the identity. We consider the preimages of $X_\sigma$ in homological grading $h(X_\sigma) - 1$. Again all these preimages will vanish since they contain a cap or a cup composed with a Jones-Wenzl projector. Thus $X_\sigma$ is not a boundary and $[X_\sigma] \neq 0$.

$\square$

**Remark 3.3.** The converse of Theorem 1.2 does not hold: One could find non zero homology classes in colored Khovanov homology whose homological grading and quantum grading does not comes from an incompressible surface as described by the statement of Theorem 1.2, see [Kin18, Section 6, Figure 10].

## 4. Example: A non-Montesinos knot

In this section we exhibit an example of a knot that is not Montesinos for which the correspondence to nonzero colored Khovanov homology classes gives an essential surface. This demonstrates the applicability of Theorem 1.1 to general classes of knots in detection of boundary slopes.

Consider the knot $K$ with diagram as shown in Figure 25, and twist region labeled $w_1, \ldots, w_5$.

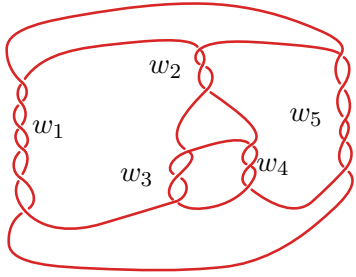

FIGURE 25. A non Montesinos knot.

Let $\mathcal{B}(K)$ be the set of boundary slopes of $K$. The *boundary slope diameter* of a knot $K$ is

$$bd(K) = max\{|s - s'| : s, s' \in \mathcal{B}(K) \setminus \{\infty\}\}.$$

We can tell that the knot is not a Montesinos knot by applying the following result from [IM08]. See [How14] for other examples of how this is used to show that a given knot is not Montesinos.

**Theorem 4.1** ([IM08]). *If $K$ is a Montesinos knot, then*

$$bd(K) \leq 2c(K),$$

*with equality if $K$ is alternating and Montesinos.*

Using SnapPy [CDGW], we see that the knot given in Figure 25 has boundary slopes that include the following rational numbers:

$$\mathcal{B}(K) \supset \{\frac{80}{1}, \frac{-79}{1}\}.$$

The diagram given has $6 + 3 \cdot 3 + 5 = 20$ crossings. Therefore $2c(K) \leq 2(20) = 40$. However, $bd(K)$ is at least $80 - (-79) = 159 > 40$ from (4), so $K$ is not a Montesinos knot.

A colored surface state on the given diagram of the knot is given in Figure 26 with slope of the corresponding normal surface equal to 4.

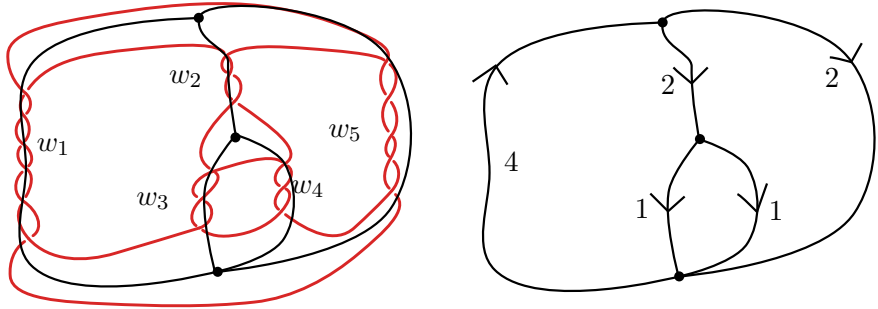

FIGURE 26. Left: The graph $G$. Right: A colored surface state.

Theorem 1.1 predicts that the slope 4 is a candidate for a boundary slope, which it is indeed as verified by SnapPy.

Department of Mathematics and Statistics, University of South Alabama, Mobile AL 36688
*Email address*: christine.rs.lee@gmail.com